\documentclass[10pt,a4paper,reqno,draft]{amsart}

\usepackage{amssymb, amsmath, textcomp, amsthm}
\usepackage{bbm,dsfont}
\usepackage[pdftex]{graphicx}
\numberwithin{equation}{section}
\title[Variation-norm Hilbert transforms along Lipschitz vector fields]{Single annulus estimates for the variation-norm Hilbert transforms along Lipschitz vector fields}
\author{Shaoming Guo}
\date{}

\usepackage[margin=1.5 in]{geometry}

\def\R{\mathbb{R}}
\def\N{\mathbb{N}}

\def\Z{\mathbb{Z}}

\def\lesim{\lesssim}

\def\mc{\mathcal}
\def\mb{\mathbbm}

\def\l{\left}
\def\r{\right}

\def\top{{\bf top}}
\def\size{{\bf size}}
\def\den{{\bf dense}}
\def\dense{\overline{{\bf dense}}}

\def\beq{\begin{equation}}
\def\endeq{\end{equation}}

\theoremstyle{plain}
\newtheorem{thm}{Theorem}[section]

\newtheorem{lem}[thm]{Lemma}

\newtheorem{defi}[thm]{Definition}

\newtheorem{rem}[thm]{Remark}

\newtheorem*{conj*}{Lipschitz Differentiation Conjecture}
\newtheorem*{openproblem*}{Open Problem}

\begin{document}
\maketitle

\begin{abstract}
Let $v$ be a planar Lipschitz vector field. We prove that the $r$-th variation-norm Hilbert transform along $v$ (defined as in \eqref{2406ee1.10}), composed with a standard Littlewood-Paley projection operator $P_k$, is bounded from $L^2$ to $L^{2, \infty}$, and from $L^p$ to itself for all $p>2$. Here $r>2$ and the operator norm is independent of $k\in \Z$. This generalises Lacey and Li's result \cite{LL1} for the case of the Hilbert transform. However, their result only assumes measurability for vector fields. In contrast to that, we need to assume vector fields to be Lipschitz.
\end{abstract}

\let\thefootnote\relax\footnote{Date: 16. Oct. 2016.  MSC (2010) 42B20, 42B25.}

\section{Introduction}

Let $v:\R^2\to S^1$ be a measurable unit vector field. Define the directional maximal operator along the vector field $v$ by
\beq
M_v f(x):=\sup_{\epsilon>0}\frac{1}{2\epsilon}\left| \int_{-\epsilon}^{\epsilon}f(x-v(x)t)dt \right|.
\endeq
Define the directional maximal operator along $v$, truncated at the scale $\epsilon_0$, by 
\beq
M_{v, \epsilon_0} f(x):=\sup_{\epsilon\le \epsilon_0}\frac{1}{2\epsilon}\left| \int_{-\epsilon}^{\epsilon}f(x-v(x)t)dt \right|.
\endeq
Similarly, we denote the directional Hilbert transform along $v$ by $H_v$, and the truncated directional Hilbert transform by $H_{v, \epsilon_0}$.\\

It is a long standing conjecture (see the discussion in Lacey and Li \cite{LL2}) that the truncated directional maximal operator and Hilbert transform along a vector field $v$ are (weakly) bounded on $L^2$ under the assumption that $v$ is Lipschitz.

In \cite{LL1}, Lacey and Li, by only assuming the vector field $v$ to be measurable but assuming the frequencies of the function $f$ to be supported on a single annulus, obtained some partial progress of the above conjecture.
\begin{thm}[\cite{LL1}]\label{2406theorem1.1}
For $k\in \Z$, let $P_k$ denote a standard $k$-th Littlewood-Paley  projection operator. For an arbitrary measurable vector field $v$, $M_v \circ P_k$ and $H_v \circ P_k$ map $L^2$ to $L^{2, \infty}$, and $L^p$ to itself for all $p>2$, with the operator norms being independent of $k\in \Z$. Moreover, both $M_v\circ P_k$ and $H_v\circ P_k$ may be unbounded on $L^p$ for any $p\le 2$.
\end{thm}

For further progress towards the above conjecture, we refer to \cite{Ba}, \cite{BT}, \cite{Bo}, \cite{CNSW}, \cite{Guo1}, \cite{LL2} and \cite{SS}.\\

In this paper, we generalise the results in Theorem \ref{2406theorem1.1} to the case of the variation-norm Hilbert transforms. Perhaps unexpectedly, the measurability assumption for the vector fields in Theorem \ref{2406theorem1.1} might not be sufficient any more. It turns out that one seems to need the conjectured optimal Lipschitz regularity. Before stating our result, we need to introduce several notations. 

For $r>2$, for a sequence of complex numbers $\{a_l: l\in \Z\}$, define its $r$-th variation-norm by
\beq
\mc{V}^r (\{a_l: l\in \Z\}):=\sup_{N\in \N; k_1<k_2<...<k_N} \l(\sum_{1\le i\le N-1}\l|a_{k_{i+1}}-a_{k_i}\r|^r\r)^{1/r}.
\endeq
Let $\psi_0:\R\to \R$ be a non-negative Schwarz function supported on $[1/2, 5/2]$ such that
\beq
\psi_0(t)=1, \forall t\in [5/4, 7/4].
\endeq 
Denote $\psi_l(t):=\psi_0(2^{-l} t)$. We could choose $\psi_0$ properly such that
\beq
\mathbbm{1}_{\R^+}(t)=\sum_{l\in \Z}\psi_l(t).
\endeq
For the vector field $v$, for $l\in \Z$, define 
\beq\label{2406ee1.8}
H_l f(x):=\int_{\R}f(x-v(x)t)\check{\psi}_l(t)dt.
\endeq
Moreover, define the variation-norm Hilbert transform along the vector field $v$ by 
\beq\label{2406ee1.9}
H^*_{v} f(x):=\mathcal{V}^r\left(\l\{\sum_{l'\le l}H_{l'} f(x): l\in \Z \r\} \right).
\endeq
Similarly, define the variation-norm Hilbert transform along the vector field $v$, truncated at the scale $\epsilon_0$, by
\beq\label{2406ee1.10}
H^*_{v, \epsilon_0} f(x):=\mathcal{V}^r\left(\l\{\sum_{-\ln \epsilon_0\le l'\le l}H_{l'} f(x): l\in \Z \r\} \right).
\endeq
Now we are ready to state our 
\begin{thm}\label{2406theorem1.2}
Let $v$ be a Lipschitz unit vector field. Let $\epsilon_0:=1/\|v\|_{Lip}$. Then $H^*_{v, \epsilon_0}\circ P_k$ maps $L^2$ to $L^{2, \infty}$ and $L^p$ to $L^p$ for any $p>2$, with the operator norms being independent of $k\in \Z$.
\end{thm}

In the above theorem, we do not know how to prove an $L^p$ bound for any $p<2$. This is supposed to be an extremely difficult problem. If one were able to prove an $L^p$ bound for $H_{v, \epsilon_0}\circ P_k$ (even the case of the Hilbert transform) with some $p<2$, then by using the almost orthogonality argument in Lacey and Li \cite{LL2} (see Chapter 5 therein), one would be able to prove the $L^2$ boundedness of $H_{v, \epsilon_0}$ by assuming $v\in C^{1+\delta}$ for some $\delta>0$.\\

The main differences between the proofs of Theorem \ref{2406theorem1.2} and Theorem \ref{2406theorem1.1} are in the tree lemma (the following Lemma \ref{1806lemma3.5}). Other than the fact that taking a variation-norm Hilbert transform destroys the orthogonality among wavelet functions, a more subtle difference is that the Lipschitz regularity will play an important role in the proof of Theorem \ref{2406theorem1.2}. As we will see, this also reflects one difference between time-frequency analysis in dimension one and higher. Moreover, the observation made in the present paper records another appearance of the conjectured Lipschitz regularity in the problem of Hilbert transforms along vector fields. \\

If we assume our vector fields to depend only on one variable, that is, $v(x_1, x_2)=(1, u(x_1))$ for some measurable function $u:\R\to \R$, then combined with Bateman's bounds on the Lipschitz-Kakeya maximal operator in \cite{Ba0} and \cite{Ba}, essential parts of the proof of Theorem \ref{2406theorem1.2} can be recycled to show
\begin{thm}\label{2909theorem1.3}
Let $u:\R\to\R$ be an arbitrary measurable function. Let $v(x_1, x_2)=(1, u(x_1))$. Then for any $p>1$, we have
\beq\label{2909ee1.9}
\|H_v^* \circ P_k f\|_p \lesim \|P_k f\|_p,
\endeq
with a constant independent of $k\in \Z$. Here $H_v^*$ is the non-truncated variation-norm Hilbert transform along $v$ defined by \eqref{2406ee1.9}.
\end{thm}

\begin{rem}
We will leave out the proof of Theorem \ref{2909theorem1.3} as the modifications needed for Bateman's argument in \cite{Ba} are almost identical to the proof of Theorem \ref{2406theorem1.2}. Lipschitz regularity is not present here due to the assumption that the vector fields are constant in the second variable $x_2$. See Remark \ref{2909remark5.5}.
\end{rem}

The bound \eqref{2909ee1.9} generalises Bateman's bound for the case of the Hilbert transforms in \cite{Ba}. Moreover, the variation-norm Hilbert transforms along vector fields appear naturally when one attempts to generalise Bateman \cite{Ba} from Hilbert transforms along vector fields to Hilbert transforms along variable polynomial curves. 

Let $i=1$ or 2. Let $u_i$ be an arbitrary measurable function. Denote $v_i(x_1, x_2)=(1, u_i(x_1))$. For a positive integer $\alpha>1$, define the Hilbert transform along the variable polynomial curve $(t, u_1(x_1)t+u_2(x_1)t^{\alpha})$ by
\beq\label{2909ee1.10}
\mathcal{H}_{\alpha} f(x_1, x_2)=\int_{\R}f(x_1-t, x_2-u_1(x_1)t-u_2(x_1)t^{\alpha})\frac{dt}{t}.
\endeq
For each fixed $k\in \Z$, by doing a partial Fourier transform in the $x_2$ variable and applying Fubini's theorem, the $L^2$ bound
\beq\label{2909ee1.11}
\|\mathcal{H}_{\alpha}\circ P_k f\|_2 \lesim \|P_k f\|_2
\endeq
is equivalent to the $L^2$ bound of the polynomial Carleson operator 
\beq
\mathcal{C}_{\alpha} f(x):=\sup_{u_1, u_2\in \R}|\int_{\R}f(x-t)e^{iu_1 t+ i u_2 t^{\alpha}}\frac{dt}{t}|,
\endeq
which has been proved by Lie \cite{Lie1} and \cite{Lie2}. \\

To prove an $L^p$ bound of the form \eqref{2909ee1.11} for certain $p$ other than 2, one natural idea is to compare the function $f(x_1-t, x_2-u_1(x_1)t-u_2(x_1)t^{\alpha})$ with $f(x_1-t, x_2-u_1(x_1)t)$ when $t$ is ``small'', by taking the advantage that the function $f$ does not oscillate fast due to its frequency localisation. When $t$ is ``large'', we expect to use certain oscillatory integral estimates similar to Stein and Wainger \cite{SteinWainger01}. By splitting $t\in \R$ into different intervals, the maximally truncated Hilbert transforms appear naturally in this problem. We will explore this idea in a future work. One similar idea has been used in the context of the Carleson-type maximal operator along curves, see \cite{GPRY}.\\

{\bf Organisation of paper:} In Section \ref{section2} we will do further reductions to the estimates in Theorem \ref{2406theorem1.2}. They are a spatial localisation and a frequency localisation separately.

In Section \ref{section3} we will introduce the time-frequency decomposition. Here we follow Section 3 in Bateman's paper \cite{Ba}.

In Section \ref{section4} we introduce several key definitions and lemmas, and in Section \ref{section5} we prove the tree lemma. This is the only place that differs from Lacey and Li's proof in \cite{LL1} (and also Bateman's proof in \cite{Ba}).\\

{\bf Notations:} Throughout this paper, we will write $x\ll y$ to mean that $x\le y/10$, $x\lesim y$ to mean that there exists a universal constant $C$ s.t. $x\le C y$, and $x\sim y$ to mean that $x\lesim y$ and $y\lesim x$.  Lastly, $\mathbbm{1}_E$ will always denote the characteristic function of the set $E$.\\

\section{Some reductions}\label{section2}

Before digging into the proof of Theorem \ref{2406theorem1.2}, we first make several reductions. Here we will take the weak type $(2, 2)$ estimate as an example. The reductions for the $L^p$ estimates can be done similarly. \\

The first reduction we will do is a spatial localisation, using the fact that we are truncating the variation-norm Hilbert transform at the scale $\epsilon_0$ in \eqref{2406ee1.10}.

By the isotropic scaling symmetry $x\to \lambda x$, we can w.l.o.g. assume that $\|v\|_{Lip}=1$. Hence $\epsilon_0=1$. For the sake of simplicity, we will use $H^*_{v}$ to stand for $H^*_{v, 1}$. What we need to prove becomes
\beq
\|H_{v}^*\circ P_{k_0}f\|_{L^{2, \infty}(\R^2)} \lesim \|P_{k_0}f\|_{L^2(\R^2)},
\endeq
with the implicit constant being independent of $k_0\in \Z$. By duality, it suffices to prove
\beq
|\langle H_{v}^*\circ P_{k_0}f, \mb{1}_E \rangle| \lesim \|P_{k_0}f\|_2 |E|^{1/2},
\endeq
for an arbitrary measurable set $E\subset \R^2$. 

For $m=(m_1, m_2)\in \Z^2$, let $\Omega_m$ denote the region 
\beq
[10^{-2}\cdot m_1, 10^{-2}\cdot (m_1+1)]\times [10^{-2}\cdot m_2, 10^{-2}\cdot (m_2+1)].
\endeq
By the triangle inequality and the Cauchy-Schwarz inequality, it suffices to prove
\beq\label{2406ee2.4}
|\langle H_{v}^*\circ P_{k_0}(P_{k_0}f\cdot \mb{1}_{\Omega_m}), \mb{1}_{E\cap \Omega_{m'}} \rangle| \lesim \frac{1}{(|m-m'|+1)^N}\|P_{k_0}f\cdot \mb{1}_{\Omega_m}\|_2 |E\cap \Omega_{m'}|^{1/2}.
\endeq
Here $N$ is certain large number. For simplicity, we will only consider the diagonal terms $m=m'$. From the forthcoming proof, it will be clear that the non-diagonal terms $m\neq m'$ can be proven similarly, and the decay in $|m-m'|$ is just from the non-stationary phase method.

Take $m=m'=0$, then \eqref{2406ee2.4} becomes
\beq\label{2406ee2.5}
|\langle H_{v}^*\circ P_{k_0}(P_{k_0}f\cdot \mb{1}_{\Omega_0}), \mb{1}_{E\cap \Omega_{0}} \rangle| \lesim \|P_{k_0}f\cdot \mb{1}_{\Omega_0}\|_2 |E\cap \Omega_{0}|^{1/2}.
\endeq
The above estimate is equivalent to 
\beq\label{2406ee2.6}
\|H_v^* P_{k_0}(P_{k_0}f\cdot \mb{1}_{\Omega_0})\|_{L^{2, \infty}(\Omega_0)}\lesim \|P_{k_0}f\cdot \mb{1}_{\Omega_m}\|_{L^2(\R^2)}.
\endeq
As we are only evaluating $H_v^* P_{k_0}$ on the domain $\Omega_0$, from now on, we can assume that our vector field $v$ is periodic in both variables, with a periodicity $(10^{-2}, 10^{-2})$. By our normalisation that $\|v\|_{Lip}=1$, we can assume that $v(x)=(1, u(x))$, with 
\beq\label{2406ee2.7}
\|u\|_{\infty}\le 10^{-2} \text{ and } \|\nabla u\|_{\infty}\le 2.
\endeq

Now we do our second reduction, which is a frequency localisation. Let $\Gamma$ denote the two-ended cone which forms an angle less than $\pi/10$ with the vertical axis. Let $\Gamma^c$ denote its complement on $\R^2$. Moreover, define
\beq
P_{\Gamma} f=\mc{F}^{-1}(\mb{1}_{\Omega}\cdot \mc{F}f).
\endeq
Here $\mc{F}$ denotes the Fourier transform, and $\mc{F}^{-1}$ its inverse. Similarly we define $P_{\Gamma^c}$. By the definition of the variation-norm Hilbert transform in \eqref{2406ee1.8} and \eqref{2406ee1.9}, and by our assumption in \eqref{2406ee2.7}, it is not difficult to see that $H_v^* P_{k_0}P_{\Gamma^c}$ is essentially the same as $P_{k_0}P_{\Gamma^c}$. Hence to prove \eqref{2406ee2.6}, what is left is to prove 
\beq
\|H_v^* P_{k_0}P_{\Gamma}(P_{k_0}f\cdot \mb{1}_{\Omega_m})\|_{L^{2, \infty}(\Omega_0)}\lesim \|P_{k_0}f\cdot \mb{1}_{\Omega_m}\|_{L^2(\R^2)}.
\endeq
Denote $P_{\Gamma_{k_0}}=P_{k_0}P_{\Gamma}$. In the following, we will focus on a slightly stronger estimate 
\beq\label{2406ee2.10}
\|H_v^* P_{\Gamma_{k_0}}f\|_{L^{2, \infty}(\R^2)}\lesim \|f\|_{L^2(\R^2)}.
\endeq
So far we have finished the steps of reductions.

\section{The time-frequency decomposition}\label{section3}
The content of this section is basically taken from Bateman's paper \cite{Ba}, with minor changes to our purpose.

For a fixed $l\in \N$, we write $\mc{D}_l$ to denote the collection of dyadic intervals of length $2^{-l}$ containded in $[-2, 2]$. Fix a smooth non-negative function $\beta:\R\to \R$ such that 
\beq
\beta(x)=1, \forall x\in [-1, 1]; \beta(x)=0, \forall |x|\ge 2.
\endeq
For each $\omega\in \mc{D}_l$, define
\beq
\beta_{\omega}(x)=\beta(2^{k_0-l}(x-c_{\omega_1})).
\endeq
Here $k_0$ is the same as in \eqref{2406ee2.10}, $\omega_1$ is the right half of $\omega$ and $c_{\omega_1}$ denotes the center of the interval $\omega_1$. Define
\beq
\beta_l(x)=\sum_{\omega\in \mathcal{D}_l} \beta_{\omega}(x).
\endeq
Note that 
\beq
\beta_l(x+2^{-l})=\beta_l(x), \forall x\in [-2, 2-2^{-l}].
\endeq
Define 
\beq
\gamma_l=\frac{1}{2}\int_{-1}^1 \beta_l(x+t)dt.
\endeq
Because of the above periodicity, we know that $\gamma_l$ is constant for $x\in [-1,1]$, independent of $l$. Say $\gamma_l(x)=\delta>0$, hence 
\beq
\frac{1}{\delta}\gamma_l(x) \mathds{1}_{[-1,1]}(x)=\mathds{1}_{[-1,1]}(x).
\endeq
Define another multiplier $\tilde{\beta}:\R\to \R$ with support in $[\frac{1}{2}, \frac{5}{2}]$ and $\tilde{\beta}(x)=1$ for $x\in [1,2]$. We define the corresponding multipliers on $\R^2$:
\begin{align*}
& \hat{m}_{\omega}(\xi, \eta)=\tilde{\beta}(2^{-k_0}\eta)\beta_{\omega}(\frac{\xi}{\eta}),\\
& \hat{m}_{l,t}(\xi, \eta)=\tilde{\beta}(2^{-k_0}\eta)\beta_l(t+\frac{\xi}{\eta}),\\
& \hat{m}_{l}(\xi, \eta)=\tilde{\beta}(2^{-k_0}\eta)\gamma_l(\frac{\xi}{\eta}).
\end{align*}

\noindent Given $l$, we would like to decompose our function $f$ accordingly:
\beq
\begin{split}
H_l\circ P_{\Gamma_{k_0}} f=(H_l\circ P_{\Gamma_{k_0}}) (m_l * f)=\int_{-1}^1 (H_l\circ P_{\Gamma_{k_0}}) (m_{l,t} * f)dt. 
\end{split}
\endeq
By the triangle inequality, 
\beq
\begin{split}
& \mathcal{V}^r \left( \l\{\int_{-1}^1\sum_{0\le l'\le l} (H_{l'}\circ P_{\Gamma_{k_0}}) (m_{l',t} * f)dt: l\in \Z\r\} \right)\\
& \lesim \int_{-1}^1 \mathcal{V}^r \left( \l\{\sum_{0\le l'\le l}(H_{l'}\circ P_{\Gamma_{k_0}}) (m_{l',t} * f): l\in \Z\r\} \right)dt.
\end{split}
\endeq
W.l.o.g. we only look at the case $t=0$. Hence, under the above notations, what we need to prove becomes 
\beq\label{1806ee2.9}
\l\| \mathcal{V}^r \left( \l\{\sum_{0\le l'\le l}\sum_{\omega\in \mc{D}_{l'}}(H_{l'}\circ P_{\Gamma_{k_0}}) (m_{\omega} * f): l\in \Z\r\} \right) \r\|_{2, \infty} \lesim \|f\|_2.
\endeq
For a fixed $\omega\in \mc{D}_{l'}$, for the frequency localised function $m_{\omega}*f$, we would also like to localise it in space. Let $\mathcal{U}_{\omega}$ be a partition of $\R^2$ by rectangles of width $2^{-k_0}$ and length $2^{-l'}$, whose long side has slope $\theta$, where $\tan \theta=-c(\omega)$. Here $c(\omega)$ denotes the center of the interval $\omega$. If $s\in \mathcal{U}_{\omega}$, we will write $\omega_s:=\omega$, and $\omega_{s,1}$ to be the right half of $\omega$, $\omega_{s, 2}$ the left half.


An element of $\mathcal{U}_{\omega}$ for some $\omega\in \mathcal{D}_l$ is called a ``tile''. Define $\varphi_{\omega}$ such that
\beq
|\hat{\varphi}_{\omega}|^2=\hat{m}_{\omega}.
\endeq
For a tile $s\in \mathcal{U}_{\omega}$, define 
\beq\label{wavelet}
\varphi_s(p):= \sqrt{|s|} \varphi_{\omega} (p-c(s)),
\endeq
where $c(s)$ is the center of $s$. Notice that 
\beq
\|\varphi_s\|_2^2= \int_{\R^2} |s| \varphi_{\omega}^2=|s|\int_{\R^2} \hat{m}_{\omega}=1,
\endeq
i.e. $\varphi_s$ is $L^2$ normalized.\\

The constructing of the tiles above by uncertainty principle is to localize the function further in space, which is realised through

\begin{lem}(See Lemma 3.1 in Page 1030 \cite{Ba})
Using notations above, we have
\beq
f* m_{\omega}(x)=\lim_{N\to \infty} \frac{1}{4N^2} \int_{[-N, N]^2} \sum_{s\in \mc{U}_{\omega}} \langle f(\cdot), \varphi_s(p+\cdot)\rangle \varphi_s(p+x)dp.
\endeq
\end{lem}
By the above lemma and the triangle inequality, to prove \eqref{1806ee2.9}, it suffices to prove 
\beq\label{1806ee2.14}
\l\|\mc{V}^r\l( \l\{\sum_{0\le l'\le l}\sum_{\omega\in \mc{D}_{l'}}\sum_{s\in \mc{U}_{\omega}}\langle f, \varphi_s\rangle \phi_s(x): l\in \Z \r\}\r)\r\|_{2, \infty} \lesim \|f\|_2.
\endeq
Here for a tile $s\in \mc{U}_{\omega}$ with $\omega\in \mc{D}_{l'}$, we have denoted
\beq
\phi_s(x):=H_{l'}\varphi_s(x).
\endeq
Now we would like to linearize the left hand side of \eqref{1806ee2.14}. Take a sequence of functions $\{a_i\}$ with 
\beq
\sum_i |a_i(x)|^{1/r^{\prime}}=1, \forall x\in \R^2.
\endeq
Take a sequence of increasing (in $i$) functions $\{k_i\}$ mapping from $\R^2$ to $\N$. Define
\beq
\mc{C} f:=\sum_{l\ge 0}\sum_{\omega\in \mathcal{D}_l}\sum_{s\in \mathcal{U}_{\omega}}\langle f, \varphi_s\rangle a_s \phi_s,
\endeq
where the coefficients $\{a_s\}$ satisfy that 
\beq\label{EE1.5}
a_s(x)=a_i(x), \text{ if } 2^{k_i(x)}\le l(s)< 2^{k_{i+1}(x)}.
\endeq
Here $l(s)$ denotes the length of the tile $s$. Hence to prove \eqref{1806ee2.14}, it suffices to prove
\beq
\|\mc{C}f\|_{2, \infty} \lesim \|f\|_2,
\endeq
with a universal constant. Moreover, by duality, we just need to prove that for an arbitrary measurable set $E\subset \R^2$, it holds that 
\beq
|\langle C f, \mathbbm{1}_E\rangle| \lesim \|f\|_2 |E|^{1/2}.
\endeq
By assuming $\|f\|_2=1$ and by using the triangle inequality, it suffices to prove
\beq\label{1806ee2.20}
\sum_{l\ge 0}\sum_{\omega\in \mathcal{D}_l}\sum_{s\in \mathcal{U}_{\omega}}|\langle f, \varphi_s\rangle \langle a_s \phi_s, \mathbbm{1}_E\rangle | \lesim  |E|^{1/2},
\endeq
which will be the main content of the rest of the current paper.

\section{Key definitions and lemmas}\label{section4}

In this section we collect the key definitions and lemmas that will be used in the proof of Theorem \ref{2406theorem1.2}.
\begin{defi}
Given two tiles $R_1, R_2\in \mc{U}$, we will write $R_1\le R_2$ if 
\beq
R_1\subset C R_2 \text{ and } \omega_{R_2}\subset \omega_{R_1}.
\endeq
Here $C$ denotes some large universal constant. 
\end{defi}

\begin{defi}
A tree is a collection $T$ of tiles with a top tile, denoted as ${\bf top}(T)$, with ${\bf top}(T)\in \mc{U}$, such that for all $s\in T$, we have $s\le \top(T)$. For $j\in \{1, 2\}$, a tree is a $j$-tree if $\omega_{\top(T)}\cap \omega_{s, j}=\emptyset$. Given a tree $T$, we will write $T_j$ to denote the maximal $j$-tree contained in $T$.
\end{defi}

For $x\in \R^2$, let $\chi$ denote the bump function 
\beq
\chi(x)=\frac{1}{1+|x|^{100}}.
\endeq
Moreover, let $\chi_s^{(p)}$ denote the $L^p$ normalised version of $\chi$ adapted to the tile $s$.

\begin{defi}\label{2406defi4.3}
For a tile $s$ and a collection of tiles $\mc{S}$, we define
\beq
E_s=\{(x, y)\in E: u(x)\in 2\cdot \omega_s\}.
\endeq
\beq
{\bf dense}(s)=\int_{E_s} \chi_s^{(1)}.
\endeq
\beq
\overline{{\bf dense}}(s)=\sup_{s'\ge s} {\bf dense}(s').
\endeq
\beq
{\bf size}(\mc{S}):=\sup_{1-\text{trees }T\subset \mc{S}}\l( \frac{1}{|{\bf top}(T)|} \sum_{s\in T}|\langle f, \varphi_s\rangle|^2 \r)^{1/2}.
\endeq
\end{defi}

\begin{rem}
Here we made a slight modification to the definition of the set $E_s$ in \cite{LL1} and \cite{Ba}. Instead of using $u\in \omega_s$, we use $u\in 2\cdot \omega_s$. This will not affect the following orthogonality lemma (Lemma \ref{1806lemma3.6}). However, it will play an important role in our \emph{Tree lemma} (Lemma \ref{1806lemma3.5}). See Remark \ref{2406remark5.6}.
\end{rem}

To prove \eqref{1806ee2.20}, by a standard argument (for example see Section 7 in \cite{Ba}), it suffices to prove the following two lemmas. 
\begin{lem}[Tree lemma]\label{1806lemma3.5}
Let $T$ be a tree. Suppose $\dense(T)\le \delta$ and $\size(T)\le \sigma$. Then 
\beq
\sum_{s\in T}|\langle f, \varphi_s\rangle \langle a_s \phi_s, \mathbbm{1}_E\rangle | \lesim \sigma \delta |\top(T)|.
\endeq
\end{lem}

\begin{lem}[Organisational lemma]\label{1806lemma3.6}
Let $\mc{S}$ be a collection of tiles. Then there exists a partition of $\mc{S}$ into trees $\mc{T}_{\delta, \sigma}$ where $\delta$ and $\sigma$ are dyadic with $\delta\lesim 1$ and $\sigma\lesim 1$ such that 
\beq\label{1806ee3.8}
\sum_{T\in \mc{T}_{\delta, \sigma}}|\top(T)|\lesim \frac{1}{\sigma^2},
\endeq
and
\beq
\sum_{T\in \mc{T}_{\delta, \sigma}}|\top(T)|\lesim \frac{|E|}{\delta}.
\endeq
\end{lem}

The proof of Lemma \ref{1806lemma3.6} remains totally the same as in \cite{Ba} and \cite{LL1}. The reason is that we are using the same definitions of size and density (up to a constant factor). Hence in the following, we will focus on the proof of Lemma \ref{1806lemma3.5}.

\section{The tree lemma}\label{section5}


In this section we will prove Lemma \ref{1806lemma3.5}. Here we follow closely the proof by Bateman \cite{Ba}. The main difference comes from the {\bf 1-tree case}. However, for the sake of completeness, we still include some of his argument. \\

First, w.l.o.g. we assume that the slope of the long side of ${\bf top}(T)$ is zero. Then, similar to the starting point of the proof of the tree lemma in Carleson's theorem in \cite{LT}, we decompose $\R^2$ into certain dyadic rectangles.

 Let $\pi_1(E)$ and $\pi_2(E)$
denote the vertical and horizontal projections of a set $E$. Let $\mc{J}_1$ be a partition of $\R$ (the horizontal axis) into dyadic intervals such that $3J\times \R$ does not contain any tile $s\in T$, and such that $J$ is maximal w.r.t. this property. In the vertical direction, we do a trivial partition, that is, we let $\mc{J}_2$ be a partition of $\R$ into intervals of width $|\pi_2({\bf top}(T))|/3$. Denote
\beq
\mc{P}=\bigcup_{J_1\in \mc{J}_1}\bigcup_{J_2\in \mc{J}_2}J_1\times J_2.
\endeq
This is a partition of $\R^2$. \\

Let us look at the term we need to bound. Note that for appropriate $\epsilon_s$ with $|\epsilon_s|=1$, we have 
\beq\label{0506ee4.3}
\begin{split}
& \sum_{s\in T} |\langle \mb{1}_F, \varphi_s\rangle\langle \mb{1}_E, a_s \phi_s\rangle| =\sum_{s\in T}\epsilon_s \langle \mb{1}_F, \varphi_s\rangle\langle \mb{1}_E, a_s \phi_s\rangle\\
&=\int_{\R^2} \sum_{s\in T} \epsilon_s \langle f, \varphi_s\rangle a_s \phi_s \mb{1}_E=\sum_{P\in \mc{P}}\int_P \sum_{s\in T} \epsilon_s \langle f, \varphi_s\rangle a_s \phi_s \mb{1}_E.
\end{split}
\endeq
Given $P\in \mc{P}$, we observe that for a tile $s\in T$, if $|\pi_1(s)|\le |\pi_1(P)|$, then by the construction of the partition $\mc{P}$, we see easily that $s\cap P=\emptyset.$ Hence when integrating over $P$, the contribution from those tiles with shorter horizontal projection is ``small''. This observation suggests the following splitting: Let 
\beq
T_P^+=\{s\in T: |\pi_1(s)|>|\pi_1(P)|\}, L_P^+=\sum_{s\in T_P^+}\epsilon_s \langle f, \varphi_s\rangle a_s \phi_s \cdot \mb{1}_E,
\endeq
and 
\beq
T_P^-=\{s\in T: |\pi_1(s)|\le |\pi_1(P)|\}, L_P^-=\sum_{s\in T_P^-}\epsilon_s \langle f, \varphi_s\rangle a_s \phi_s \cdot \mb{1}_E.
\endeq
Then
\beq\label{0506ee4.6}
\eqref{0506ee4.3}= \sum_{P\in \mc{P}} \int_P L_P^- + \sum_{P\in \mc{P}} \int_P L_P^+.
\endeq
As we have explained above, the former term on the right hand side of the last expression has ``small'' contribution. It can be bounded exactly in the same way as in Subsection 11.1 of Bateman \cite{Ba}. Hence we leave out the details.

Now we turn to the latter term in \eqref{0506ee4.6}. By measuring how far $P$ is from ${\bf top}(T)$, we do a further splitting. For $k\ge 1$, define 
\beq\label{0506ee4.7}
\begin{split}
& \mc{P}_0 =\{P\in \mc{P}: \frac{\text{dist}(\pi_2(P), \pi_2({\bf top}(T)))}{|\pi_2({\bf top}(T))|}\le 1\},\\
& \mc{P}_k =\{P\in \mc{P}: \frac{\text{dist}(\pi_2(P), \pi_2({\bf top}(T)))}{|\pi_2({\bf top}(T))|}\in (2^{k-1}, 2^k]\}.
\end{split}
\endeq

\noindent We show that the term $L^+_P$ has small support. Precisely speaking,
\begin{lem}\label{0506lemma4.2}(Claim 11.1 in \cite{Ba})
For $P\in \mc{P}_k$, if we use $E_p$ to denote the support of the function $L_P^+\cdot \mb{1}_P$, then $|E_P|\lesim 2^{100k}\delta |P|$.
\end{lem}
\noindent {\bf Proof of Lemma \ref{0506lemma4.2}:} Here we give a direct proof instead of arguing by contraction as in \cite{Ba}. By the construction of $\mc{P}$, we know that there is some $s\in T$ with $|s|\sim |P|$ such that $s\subset C 2^k P$ for some universal constant $C$. Moreover, among these tiles, we let $s_0$ denote the one with the smallest area. Hence 
\beq
\forall s\in T_P^+: \omega_s\subset \omega_{s_0}.
\endeq
This further implies that
\beq
\text{supp} (L_P^+\cdot \mb{1}_P)\subset P\cap E_{s_0}.
\endeq
By the assumption that ${\bf dense}(s_0)\lesim \delta$, we obtain
\beq
\delta\gtrsim \int_{E_{s_0}} \chi_s^{(1)}\ge \int_{P\cap E_{s_0}} \chi_s^{(1)},
\endeq
which further implies that 
\beq
|P\cap E_{s_0}|\lesim 2^{100k}\delta |P|.
\endeq
So far we have finished the proof of Lemma \ref{0506lemma4.2}. $\Box$\\

We proceed with the estimate of the second term of \eqref{0506ee4.6}, which is 
\beq
\sum_{P\in \mc{P}}\int_P L_P^+\cdot \mb{1}_{E_P}.
\endeq
Clearly for every $s\in T$, either $\omega_{s, 1}\cap \omega_{\top(T)}=\emptyset$ or $\omega_{s, 2}\cap \omega_{\top(T)}=\emptyset$, so our tree $T$ can be partitioned as $T=T_1\cup T_2$, where $T_j$ is a $j$-tree. Let 
\beq
(T_P^+)_j= T_P^+\cap T_j.
\endeq
Notice that $(T_P^+)_j$ is still a $j$-tree.\\

{\bf The $2$-tree case:} This case is relatively easier to handle, due to the separation of the supports of $\{\phi_s\}_{s\in T_2}$. Precisely speaking, for $s, t\in T_2$, if there exists a point $x\in \R^2$ such that $\phi_s(x)\phi_t(x)\neq 0$, then $|s|=|t|$. The rest of the argument for this case remains the same as in Subsection 11.2.1 of Bateman \cite{Ba}, and again we leave out the details. 

{\bf The $1$-tree case:} This case is where the main difference between the Hilbert transform and the variation-norm Hilbert transform lies. In this case we need to appeal to the orthogonality of the wavelet functions $\{\phi_s\}_{s\in T_1}$. When applying a variation-norm, this orthogonality will be destroyed. To recover such an orthogonality, one straightforward idea is to use the boundedness of the variation-norm Hilbert transform in dimension one. However, this turns out to be not enough. The Lipschitz assumption on the vector field will start to play an important role. \\

Let
\beq\label{0606ee4.26}
\alpha_s(x)=\int_{\R} \psi_s(t) \varphi_s(x_1-t, x_2)dt.
\endeq
In the following, we would like to compare $\alpha_s$ with $\phi_s$ at those points $x\in P$ such that $\phi_s\neq 0$. Recall that 
\beq\label{0606ee4.27}
\phi_s(x)=\int_{\R} \psi_s(t)\varphi_s(x_1-t, x_2-t u(x_1))dt.
\endeq
In the definition of $\alpha_s$, we are integrating along the direction $(1, 0)$. This is determined by our assumption at the beginning of the proof of the tree lemma that the long side of $\top(T)$ is parallel to $(1, 0)$. If the function $\phi_s$ does not vanish at a point $x$, then necessarily $u(x_1)$ is ``small'', or more precisely,  
\beq
|u(x_1)|\le 2^{-k_0}/|\pi_1(s)|\sim w(s)/l(s).
\endeq
Hence one might expect that by the fundamental theorem of calculus, $|\alpha_s-\phi_s|$ could be bounded by some ``small'' factor.\\

We proceed with the details. Following the above idea, we write 
\beq
\begin{split}
& \int_P \sum_{s\in (T_P^+)_1} \epsilon_s \langle f, \varphi_s\rangle a_s\phi_s\cdot \mb{1}_{E_P}\\
& =\int_P \sum_{s\in (T_P^+)_1} \epsilon_s \langle f, \varphi_s\rangle a_s\alpha_s\cdot \mb{1}_{E_P} + \int_P \sum_{s\in (T_P^+)_1} \epsilon_s \langle f, \varphi_s\rangle a_s(\phi_s-\alpha_s)\cdot \mb{1}_{E_P},
\end{split} 
\endeq
which will be called $I_P$ and $II_P$ separately. \\

The estimate of the term $II_P$ is slightly easier. We need
\begin{lem}\label{0606lemma4.3}(See Claim 11.3 in Page 1052 \cite{Ba})
For $P\in \mc{P}_k$, for all $x\in P$, we have the pointwise estimate 
\beq
|\sum_{s\in (T_P^+)_1} \epsilon_s \langle f, \varphi_s\rangle a_s(\phi_s-\alpha_s)(x)\cdot \mb{1}_{E_P} (x)| \lesim 2^{-200k}\sigma.
\endeq
\end{lem}
\noindent The proof of the above lemma is done simply by using the expressions \eqref{0606ee4.26} and \eqref{0606ee4.27}, and applying the fundamental theorem of calculus to the second variable of the function $\varphi_s$. Lemma \ref{0606lemma4.3}, combined with Lemma \ref{0506lemma4.2}, gives 
\beq
\sum_k \sum_{P\in \mc{P}_k} II_P \lesim \sum_k \sum_{P\in \mc{P}_k} 2^{-100k} \sigma \delta |P|\lesim \sigma \delta |\top(T)|.
\endeq
This finishes the estimate of the term $II_P$.\\

What remains is to estimate the term $I_P$. By the definition of the variation-norm, we observe that 
\beq\label{2206ee4.31}
\l|\sum_{s\in (T_P^+)_1} \epsilon_s \langle f, \varphi_s\rangle a_s\alpha_s\r| \le \mc{V}^r\l(\l\{\sum_{s\in (T_P^+)_1,1\le |\pi_1(s)|\le 2^l} \epsilon_s \langle f, \varphi_s\rangle \alpha_s: l\in \Z\r\}\r).
\endeq
We need to integrate the right hand side of the above expression over $P$. Here comes the main difference with the proof of the tree lemma in \cite{Ba} and \cite{LL1}. In the case of the Hilbert transforms along vector fields, one bounds the left hand side of \eqref{2206ee4.31} by a strong maximal function. Moreover, this maximal function is essentially a constant on $P$. 

In the case of the variation-norm Hilbert transforms, the above idea is no longer sufficient. The right hand side of \eqref{2206ee4.31} does not behave like a constant on the whole $P$ any more. It turns out that only on every horizontal line should one view the right hand side of \eqref{2206ee4.31} as a constant. 

We proceed with the details. For each $y\in \pi_2(P)$, let $l_y$ denote the horizontal line that passes through $(0, y)$. Write 
\beq
\begin{split}
& \int_P \mc{V}^r\l(\l\{\sum_{s\in (T_P^+)_1, |\pi_1(s)|\le 2^l} \epsilon_s \langle f, \varphi_s\rangle \alpha_s: l\in \Z\r\}\r)\\
&=\int_{\pi_2(P)}\int_{P\cap l_y} \mc{V}^r\l(\l\{\sum_{s\in (T_P^+)_1, |\pi_1(s)|\le 2^l} \epsilon_s \langle f, \varphi_s\rangle \alpha_s(x, y): l\in \Z\r\}\r)dxdy.
\end{split}
\endeq
We fix a $y_0\in \pi_2(P)$ and look at the integration over $P\cap l_{y_0}$. Take $x_0\in \pi_1(P\cap l_{y_0})$. At the point $(x_0, y_0)$, suppose that the supremum in the definition of the variatoin-norm
\beq
\mc{V}^r \l(\l\{\sum_{s\in (T_P^+)_1, |\pi_1(s)|\le 2^l} \epsilon_s \langle f, \varphi_s\rangle \alpha_s(x_0, y_0): l\in \Z\r\}\r)
\endeq 
is obtained in the form of 
\beq
\sum_{s\in (T_P^+)_1}\epsilon_s \langle f, \varphi_s\rangle a_s(x_0, y_0)\alpha_s(x_0, y_0),
\endeq
where $\{a_s\}_{s\in (T_P^+)_1}$ is such that, there exists an increasing sequence $\{\tilde{k}_i\}_{i\in \N}\subset \Z$ such that, 
\beq
\text{if } 2^{\tilde{k}_i}\le |\pi_1(s)|<2^{\tilde{k}_{i+1}}, \text{ then }a_s(x_0, y_0)=a_i(x_0, y_0), \text{ with } \sum_i |a_i(x_0, y_0)|^r=1.
\endeq
When calculating the $\mc{V}^r$ norm at another point $(x, y_0)\in P$, we use the fact that the function $\alpha_s$ behaves like a constant on $P\cap l_{y_0}$. To be precise, suppose that there exists a sequence $\{k_i\}_{i\in \N}$ such that 
\beq
\begin{split}
& \mc{V}^r \l(\l\{\sum_{s\in (T_P^+)_1, |\pi_1(s)|\le 2^l} \epsilon_s \langle f, \varphi_s\rangle \alpha_s(x, y_0): l\in \Z\r\}\r)\\
&=\l( \sum_{i\in \N} \l( \sum_{s: 2^{k_i}\le |\pi_1(s)|<2^{k_{i+1}}} \epsilon_s \langle f, \varphi_s\rangle \alpha_s(x, y_0) \r)^r \r)^{1/r}.
\end{split}
\endeq
By the triangle inequality, the right hand side of the above expression can be bounded by
\beq\label{2206ee4.37}
\begin{split}
&\l( \sum_{i\in \N} \l( \sum_{s: 2^{k_i}\le |\pi_1(s)|<2^{k_{i+1}}} \epsilon_s \langle f, \varphi_s\rangle (\alpha_s(x, y_0)-\alpha_s(x_0, y_0)) \r)^r \r)^{1/r} \\
& + \l( \sum_{i\in \N} \l( \sum_{s: 2^{k_i}\le |\pi_1(s)|<2^{k_{i+1}}} \epsilon_s \langle f, \varphi_s\rangle \alpha_s(x_0, y_0) \r)^r \r)^{1/r},
\end{split}
\endeq
which will be denoted by $III$ and $IV$ separately. We need to integrate the above expression over $P$. Concerning the term $III$, we have
\beq
\int_{\pi_2(P)}\int_{P\cap l_y} III\cdot \mb{1}_{E_P} \lesim \int_{\pi_2(P)}\int_{P\cap l_y} \sum_{s\in (T_P^+)_1}|\langle f, \varphi_s\rangle||\alpha_s(x, y)-\alpha_s(x_0, y)|\cdot \mb{1}_{E_P}.
\endeq
For each $s\in (T_P^+)_1$, by the fundamental theorem in calculus, we obtain
\beq
|\alpha_s(x, y)-\alpha_s(x_0, y)|\lesim 2^{-200k} \frac{1}{\sqrt{|s|}} \cdot \frac{|\pi_1(P)|}{|\pi_1(s)|}.
\endeq
This, together with the bound $|\langle f, \varphi_s\rangle|\lesim \sigma \sqrt{|s|}$ and the estimate for $|E_P|$ in Lemma \ref{0506lemma4.2}, implies that 
\beq
\int_{\pi_2(P)}\int_{P\cap l_y} III\cdot \mb{1}_{E_P} \lesim 2^{-100k}\sigma \delta |P|.
\endeq
By summing over $k$ and $P$, we obtain the desired upper bound $\sigma \delta |\top(T)|$. Hence we have finished the proof of the term $III$.\\

We turn to the latter term in \eqref{2206ee4.37}. 
\beq\label{2306ee4.41}
\begin{split}
&\int_{\pi_2(P)}\int_{P\cap l_y} IV\cdot \mb{1}_{E_P}\\
&=\int_{\pi_2(P)}\int_{P\cap l_y}  \l( \sum_{i\in \N} \l( \sum_{s: 2^{k_i}\le |\pi_1(s)|<2^{k_{i+1}}} \epsilon_s \langle f, \varphi_s\rangle \alpha_s(x_0, y) \r)^r \r)^{1/r}\cdot \mb{1}_{E_P}.
\end{split}
\endeq
To proceed, we need
\begin{lem}\label{2406lemma5.3}
Under the above notations, we have that 
\beq\label{2406ee5.42}
|E_P\cap l_y|\lesim 2^{100k} \delta |P\cap l_y|,
\endeq
for each $y\in \pi_2(P)$.
\end{lem}
\begin{rem}
The above lemma strengthens the estimate for $|E_P|$ in Lemma \ref{0506lemma4.2}. The estimate \eqref{2406ee5.42} says that on each horizontal line segment $P\cap l_y$, the support of $L_P^+$, occupies at most a proportion of $2^{100k}\delta$. Apparently, this can not be true for arbitrary measurable vector field. Indeed, as we will see in the following proof, the Lipschitz regularity will play a crucial role.
\end{rem}
\begin{rem}\label{2909remark5.5}
If we assume that our vector field is constant along vertical lines, then it is easy to see that Lemma \ref{2406lemma5.3} is equivalent with Lemma \ref{0506lemma4.2}.
\end{rem}
We postpone the proof of Lemma \ref{2406lemma5.3} to the end of this subsection, and continue with the estimate on the term \eqref{2306ee4.41}. By Lemma \ref{2406lemma5.3}, the right hand side of \eqref{2306ee4.41} can be bounded by 
\beq
 \int_{\pi_2(P)} 2^{100k}\delta \int_{P\cap l_y} \l( \sum_{i\in \N} \l( \sum_{s: 2^{k_i}\le |\pi_1(s)|<2^{k_{i+1}}} \epsilon_s \langle f, \varphi_s\rangle \alpha_s(x_0, y) \r)^r \r)^{1/r},
\endeq
which can be further bounded by
\beq\label{2306ee4.42}
\begin{split}
& \lesim \int_{\pi_2(P)}2^{100k} \delta \int_{P\cap l_y} \l( \sum_{i\in \N} \l( \sum_{s: 2^{k_i}\le |\pi_1(s)|<2^{k_{i+1}}} \epsilon_s \langle f, \varphi_s\rangle \alpha_s(x, y) \r)^r \r)^{1/r}\\
& + \int_{\pi_2(P)} 2^{100k}\delta \int_{P\cap l_y} \l( \sum_{i\in \N} \l( \sum_{s: 2^{k_i}\le |\pi_1(s)|<2^{k_{i+1}}} \epsilon_s \langle f, \varphi_s\rangle (\alpha_s(x, y)-\alpha_s(x_0, y)) \r)^r \r)^{1/r}.
\end{split}
\endeq
The latter term can be bounded in the same way as the former term $III$ in \eqref{2206ee4.37}, hence we leave out the details. By the definition of the variation-norm, the former term in \eqref{2306ee4.42} can be bounded by
\beq\label{2306ee4.43}
\begin{split}
& 2^{100k} \delta \int_P \mc{V}^r\l(\l\{\sum_{s\in (T_P^+)_1, |\pi_1(s)|\le 2^l} \epsilon_s \langle f, \varphi_s\rangle \alpha_s: l\in \Z\r\}\r)\\
& \lesim 2^{100k} \delta \int_P \mc{V}^r\l(\l\{\sum_{s\in T, |\pi_1(s)|\le 2^l} \epsilon_s \langle f, \varphi_s\rangle \alpha_s: l\in \Z\r\}\r).
\end{split}
\endeq

\noindent Notice that in the last term of the above expression, the function that we are integrating over $P$ does not depend on $P$ any more. For simplicity, we denote
\beq
F_{T, l}:=\sum_{s\in T, |\pi_1(s)|\le 2^l} \epsilon_s \langle f, \varphi_s\rangle \alpha_s.
\endeq
Hence the last term in \eqref{2306ee4.43} becomes 
\beq\label{0606ff4.35}
2^{100k}\delta \int_P \mc{V}^r\l(\l\{F_{T, l}: l\in \Z \r\} \r).
\endeq
For a fixed $k\in \N$, we sum over $P\in \mc{P}_k$ to obtain
\beq\label{0606ee4.34}
\eqref{0606ff4.35}\lesim  2^{100k} \delta \int_{\bigcup_{P\in\mc{P}_k}P} \mc{V}^r\l(\l\{F_{T, l}: l\in \Z\r\}\r).
\endeq
Notice that $|\bigcup_{P\in \mc{P}_k} P|\sim |\top(T)|$, hence applying H\"older's inequality to the right hand side of \eqref{0606ee4.34}, we further obtain
\beq
\eqref{2306ee4.43}\lesim 2^{100k}\delta |\top(T)|^{1/2} \l(\int_{\bigcup_{P\in\mc{P}_k}P} \l|\mc{V}^r\l(\l\{F_{T, l}: l\in \Z\r\}\r)\r|^2\r)^{1/2}
\endeq
Hence what remains is to prove
\beq
\l(\int_{\bigcup_{P\in\mc{P}_k}P} \l|V^r\l(\l\{F_{T, l}: l\in \Z\r\}\r)\r|^2\r)^{1/2} \lesim 2^{-Nk}\sigma |\top(T)|^{1/2}.
\endeq
However, this can be proven simply by using the variation-norm estimate for the one-dimensional Hilbert transform (see for example \cite{CJRW} and \cite{JSW}) and Fubini's theorem. The decaying factor $2^{-Nk}$ can again be obtained by the non-stationary phase method. We leave out the details.\\

{\bf Proof of Lemma \ref{2406lemma5.3}:} We argue by contradiction. In the proof of Lemma \ref{0506lemma4.2}, we have found a tile $s_0\in T$ with 
\beq\label{2406ee5.51}
|s_0|\sim |P|, s_0\subset 2^k P \text{ and } E_P\subset P\cap E_{s_0}.
\endeq
Suppose that for any large number $M$, there exists $y_0\in \pi_2(P)$ such that 
\beq
|E_P\cap l_{y_0}|\ge M 2^{100k}\delta |P\cap l_y|.
\endeq
Then \eqref{2406ee5.51} implies that 
\beq
|P\cap E_{s_0}\cap l_{y_0}|\ge M 2^{100k}\delta |P\cap l_y|.
\endeq
For a point $(x_0, y_0)\in P\cap E_{s_0}\cap l_{y_0}$, by the assumption that $\|v\|_{Lip}=1$, it is easy to see that 
\beq\label{2406ee5.54}
u(x_0, y)\in 2\cdot \omega_{s_0}, \forall (x_0, y)\in P.
\endeq
\begin{rem}\label{2406remark5.6}
The above estimate will not be true if we do not modify the definition of the set $E_{s_0}$ in Definition \ref{2406defi4.3}.
\end{rem}
The estimate \eqref{2406ee5.54}, together with the definition of $\den(s_0)$ in Definition \ref{2406defi4.3}, implies that 
\beq
\den(s)\ge M \delta.
\endeq
This is a contradiction to our assumption that $\dense(T)\lesim \delta$. Hence we have finished the proof of Lemma \ref{2406lemma5.3}, and thus the proof of Theorem \ref{2406theorem1.2}. $\Box$

\noindent Shaoming Guo, Institute of Mathematics, University of Bonn\\
Email: shaoguo@iu.edu\\
Address: Endenicher Allee 60, 53115, Bonn, Germany\\
Current Address: 831 E Third St, Bloomington, IN 47405.

\end{document}